\documentclass{amsart}
\usepackage{lastpage}
\usepackage{yhmath,verbatim,mathabx}
\usepackage[all]{xy}

\DeclareFontFamily{U}{mathx}{\hyphenchar\font45}
\DeclareFontShape{U}{mathx}{m}{n}{
	<5> <6> <7> <8> <9> <10>
	<10.95> <12> <14.4> <17.28> <20.74> <24.88>
	mathx10
}{}
\DeclareSymbolFont{mathx}{U}{mathx}{m}{n}
\DeclareFontSubstitution{U}{mathx}{m}{n}
\DeclareMathAccent{\widecheck}{0}{mathx}{"71}
\DeclareMathAccent{\wideparen}{0}{mathx}{"75}

\newcommand{\bdis}{\begin{displaymath}}
\newcommand{\edis}{\end{displaymath}}
\newcommand{\be}{\begin{equation}}
\newcommand{\ee}{\end{equation}}
\newcommand{\mbb}{\mathbb}
\newcommand{\mcal}{\mathcal}

\newcommand{\zf}{\zeta\left(\frac{1}{2}+it\right)}

\theoremstyle{definition}

\theoremstyle{remark}
\newtheorem{remark}[]{Remark}

\newtheorem*{mydef1}{{\bf Theorem}}

\newtheorem*{mydef41}{{\bf Corollary 1}}

\newtheorem*{mydef42}{{\bf Corollary 2}}

\newtheorem*{mydef43}{{\bf Corollary 3}}

\newtheorem*{mydef5}{{\bf Lemma}}

\numberwithin{equation}{section}

\begin{document}

\title{Jacob's ladders and the set of elementary meta-functional equations giving new kind of interactions of $\zeta(s)$ with itself} 

\author{Jan Moser}

\address{Department of Mathematical Analysis and Numerical Mathematics, Comenius University, Mlynska Dolina M105, 842 48 Bratislava, SLOVAKIA}

\email{jan.mozer@fmph.uniba.sk}


\begin{abstract}
In this paper we obtain the set of grafts from some set of 12 elementary functions on 12 critical strips. We make use this directly for grafting of elements of the corresponding set of $\zeta$-factorization formulas. This procedure gives the set of elementary meta-functional equations that represents the set of elementary interactions of the Riemann's zeta-function with itself. 

\end{abstract}
\maketitle

\section{Introduction} 

\subsection{} 

Let us recall that the proof of meta-functional equation given in \cite{9} is based: 
\begin{itemize}
	\item[(a)] on new notions and methods in the theory of Riemann's zeta-function we have introduced in our series of 50 papers concerning Jacob's ladders, these can be found an arXiv [math.CA] starting with the paper \cite{1}, 
	\item[(b)] on the classical H. Bohr's theorem in analysis (1914) concerning infinite set of roots of equation 
	\be \label{1.1} 
	\zeta(s)=a,\ a\not=0,\ a\in\mbb{C},\ \frac 12<\sigma_1<\sigma<\sigma_2<1;\ s=\sigma+it, 
	\ee 
	\item[(c)] on our notion of grafting of complete hybrid formula (see \cite{9}) that serves for synthesis of conceptions (a) and (b). 
\end{itemize} 

\subsection{} 

In this paper we give another type of the proof. Namely, we make the following operations on a given set of elementary functions: 
\begin{itemize}
	\item[(a)] we construct corresponding set of $\zeta$-factorization formulas according to our algorithm from \cite{4} and \cite{5}, 
	\item[(b)] we make grafting (see \cite{9}) of every element of the set (a). 
\end{itemize}

\begin{remark}
Corresponding set of elementary meta-functional equations is obtained in result of above mentioned operations. All these formulas are, simultaneously, synergetic ones, comp. \cite{8}. 
\end{remark}

\subsection{} 

We have obtained, for example: 
\begin{itemize}
	\item[(a)] elementary meta-functional equation 
	\be \label{1.2} 
	\begin{split}
	& |\zeta(w_6^n)|\prod_{r=1}^{k_6}\frac{|\zeta(\frac 12+i\alpha_r^{6,k_6,n})|^2}{|\zeta(\frac 12+i\beta_r^{k_6,n})|^2}\sim \\ 
	& \frac{5}{16}+\frac{15}{32}|\zeta(w_{10}^n)|+\frac{3}{16}|\zeta(w_{11}^n)|+\frac{1}{32}|\zeta(w_{12}^n)|,\ L\to\infty, 
	\end{split}
	\ee 
	\item[(b)] more complicated meta-functional equation 
	\be \label{1.3} 
	\begin{split}
	& 2\left\{
	|\zeta(w_6^n)|\prod_{r=1}^{k_6}\frac{|\zeta(\frac 12+i\alpha_r^{6,k_6,n})|^2}{|\zeta(\frac 12+i\beta_r^{k_6,n})|^2}+
	|\zeta(w_7^n)|\prod_{r=1}^{k_7}\frac{|\zeta(\frac 12+i\alpha_r^{7,k_7,n})|^2}{|\zeta(\frac 12+i\beta_r^{k_7,n})|^2}
	\right\}+1\sim \\ 
	& 3\left\{
	|\zeta(w_3^n)|\prod_{r=1}^{k_3}\frac{|\zeta(\frac 12+i\alpha_r^{3,k_3,n})|^2}{|\zeta(\frac 12+i\beta_r^{k_3,n})|^2}+
	|\zeta(w_4^n)|\prod_{r=1}^{k_4}\frac{|\zeta(\frac 12+i\alpha_r^{4,k_4,n})|^2}{|\zeta(\frac 12+i\beta_r^{k_4,n})|^2}
	\right\},\ L\to\infty 
	\end{split}
	\ee 
	with notations according to those used in \cite{9}. 
\end{itemize} 

\begin{remark}
Formulae (\ref{1.2}) and (\ref{1.3}) are asymptotic forms of corresponding exact formulae. 
\end{remark} 

\begin{remark}
Let us note that new type result concerning interaction of the Riemann's zeta-functions with itself  is expressed by any of two synergetic formulae (\ref{1.2}) and (\ref{1.3}). 
\end{remark} 

\section{Set of $\zeta$-factorization formulas} 

\subsection{} 

In this paper we use the following two sets of elementary functions 
\be \label{2.1} 
\begin{split}
& \{f_1(t),\dots,f_9(t)\}=\\ 
& =\{\sin^2t,\cos^2t,\sin^4t,\cos^4t,\sin^6t,\cos^6t,\cos 2t,\cos4t,\cos6t\},\\ 
& t\in [\pi L,\pi L+U],\ L\in\mbb{N}, 
\end{split}
\ee  
and 
\be \label{2.2} 
\begin{split}
& \{f_{10}(t),f_{11}(t),f_{12}(t)\}=
\left\{
\frac{\sin 2U}{2U},\frac{\sin 4U}{4U}, \frac{\sin 6U}{6U}
\right\}, \\ 
& U\in (0,\pi/12). 
\end{split}
\ee 
By making use our algorithm of generating $\zeta$-factorization formulas (see \cite{5}, (3.1) -- (3.11), comp. \cite{4}) on the set (\ref{2.1}) we obtain the following. 

\begin{mydef5}
\be \label{2.3} 
\sin^2\alpha_0^{1,k_1}\prod_{r=1}^{k_1}\frac{\tilde{Z}^2(\alpha_r^{1,k_1})}{\tilde{Z}^2(\beta_r^{k_1})} = \frac 12-\frac 12\frac{\sin 2U}{2U}, 
\ee 

\be \label{2.4} 
\cos^2\alpha_0^{2,k_2}\prod_{r=1}^{k_2}\frac{\tilde{Z}^2(\alpha_r^{2,k_2})}{\tilde{Z}^2(\beta_r^{k_2})} = \frac 12+\frac 12\frac{\sin 2U}{2U}, 
\ee 

\be \label{2.5} 
\sin^4\alpha_0^{3,k_3}\prod_{r=1}^{k_3}\frac{\tilde{Z}^2(\alpha_r^{3,k_3})}{\tilde{Z}^2(\beta_r^{k_3})} = \frac 38-\frac{\sin 2U}{4U}+\frac{\sin 4U}{32U}, 
\ee 

\be \label{2.6} 
\cos^4\alpha_0^{4,k_4}\prod_{r=1}^{k_4}\frac{\tilde{Z}^2(\alpha_r^{4,k_4})}{\tilde{Z}^2(\beta_r^{k_4})} = \frac 38+\frac{\sin 2U}{4U}+\frac{\sin 4U}{32U}, 
\ee 

\be \label{2.7} 
\sin^6\alpha_0^{5,k_5}\prod_{r=1}^{k_5}\frac{\tilde{Z}^2(\alpha_r^{5,k_5})}{\tilde{Z}^2(\beta_r^{k_5})} = \frac{5}{16}-\frac{15}{64}\frac{\sin 2U}{U}+\frac{3}{64}\frac{\sin 4U}{U}-\frac{1}{192}\frac{\sin 6U}{U}, 
\ee 

\be \label{2.8} 
\cos^6\alpha_0^{6,k_6}\prod_{r=1}^{k_6}\frac{\tilde{Z}^2(\alpha_r^{6,k_6})}{\tilde{Z}^2(\beta_r^{k_6})} = \frac{5}{16}+\frac{15}{64}\frac{\sin 2U}{U}+\frac{3}{64}\frac{\sin 4U}{U}+\frac{1}{192}\frac{\sin 6U}{U}, 
\ee 

\be \label{2.9} 
\cos(2\alpha_0^{7,k_7})\prod_{r=1}^{k_7}\frac{\tilde{Z}^2(\alpha_r^{7,k_7})}{\tilde{Z}^2(\beta_r^{k_7})} = \frac{\sin 2U}{2U}, 
\ee 

\be \label{2.10} 
\cos(4\alpha_0^{8,k_8})\prod_{r=1}^{k_8}\frac{\tilde{Z}^2(\alpha_r^{8,k_8})}{\tilde{Z}^2(\beta_r^{k_8})} = \frac{\sin 4U}{4U}, 
\ee 

\be \label{2.11} 
\cos(6\alpha_0^{9,k_9})\prod_{r=1}^{k_9}\frac{\tilde{Z}^2(\alpha_r^{9,k_9})}{\tilde{Z}^2(\beta_r^{k_9})} = \frac{\sin 6U}{6U}, 
\ee 
\bdis 
\forall\- L\geq L_0>0,\ 1\leq k_1,\dots,k_9\leq k_0, 
\edis 
(here we fix arbitrary $k_0\in\mbb{N}$), where 
\be \label{2.12} 
\begin{split}
& \alpha_r^{l,k_l}=\alpha_r(U,\pi L,k_l,f_l),\ r=0,1,\dots,k_l,\ l=1,\dots,9, \\ 
& \beta_r^{k_l}=\beta_r(U,\pi L,k_l),\ r=1,\dots,k_l, \\ 
& \alpha_0^{l,k_l}\in (\pi L,\pi L+U),\ 
\alpha_r^{l,k_l},\beta_r^{k_l}\in 
(\overset{r}{\wideparen{\pi L}},\overset{r}{\wideparen{\pi L+U}}), \\ 
& r=1,\dots,k_l, 
\end{split}
\ee  
and 
\bdis 
(\overset{r}{\wideparen{\pi L}},\overset{r}{\wideparen{\pi L+U}})
\edis  
is the $r$-th reverse iteration of the interval $(\pi L,\pi L+U)$ by means of the Jacob's ladder (see \cite{3}) and 
\be \label{2.13} 
\tilde{Z}^2(t)=\frac{|\zf|^2}{\omega(t)},] \omega(t)=\left\{1+\mcal{O}\left(\frac{\ln\ln t}{\ln t}\right)\right\}\ln t, 
\ee 
(see \cite{2}, (6.1), (6.7), (7.7), (7.8) and (9.1)). 
\end{mydef5} 

\section{The set of corresponding grafts} 

\subsection{} 

We insert into the basic strip 
\bdis 
(\sigma_1,\sigma_2)\times (0,+\infty),\ \frac 12<\sigma_1<\sigma_2<1 
\edis 
the following twelve strips 
\be \label{3.1} 
S^l_{\sigma_0^l,\delta}=(\sigma_0^l-\delta,\sigma_0^l+\delta)\times (0,+\infty),\ l=1,\dots,12, 
\ee 
where 
\be \label{3.2} 
\frac{1}{2}<\sigma_1<\sigma_0^1<\sigma_0^2<\dots<\sigma_0^{12}<\sigma_2<1, 
\ee 
and 
\be \label{3.3} 
\sigma_1<\sigma_0^1-\delta,\ \sigma_0^{12}+\delta<\sigma_2,\ \sigma_0^l+\delta<\sigma_0^{l+1}-\delta,\ l=1,\dots,11  
\ee 
for suffciently small $\delta>0$. Of course, 
\be \label{3.4} 
S^k_{\sigma_0^k,\delta}\bigcap S^l_{\sigma_0^l,\delta}=\emptyset,\ k\not=l,\ k,l=1,\dots,12. 
\ee  

\subsection{} 

Next, we choose arbitrary finite set (comp. (\ref{2.2})) 
\be \label{3.5} 
\{U_n\}:\ 0<U_1<U_2<\dots<U_{n_0}<\frac{\pi}{12} 
\ee 
that fulfills the condition (for example) 
\be \label{3.6} 
U_{n+1}-U_n>10^{-34},\ n=1,\dots,n_0-1,\ U_1,\frac{\pi}{12}-U_{n_0}>10^{-43}. 
\ee 

\begin{remark}
Condition (\ref{3.6}) corresponds with the point of view of Jakov Zeldovich for using maths to study real-world phenomena. Namely, the choice done in (\ref{3.6}) follows from the fact that Planck's length and Planck's time have the following values 
\bdis 
L_p=8.1\times 10^{-35}\mbox{cm},\ T_p=2.7\times 10^{-43}\mbox{s}. 
\edis 
\end{remark}

\begin{remark}
Of course, 
\be \label{3.7} 
\frac{\pi}{12}:10^{-43}<10^{43} \ \Rightarrow\ n_0=n_0(L)\leq 10^{43},\ L\to\infty, 
\ee 
(see (\ref{3.5}), (\ref{3.6})). 
\end{remark} 

\subsection{} 

Let $\{U_n\}$ be some admissible set. We use the following notations (comp. (\ref{2.12})) 
\be \label{3.8} 
\begin{split}
& \alpha_0^{l,k_l,n}=\alpha_0(U_n,\pi L,k_l,f_l),\ l=1,\dots,9, \\ 
& \vdots 
\end{split}
\ee  
Since $U_n\in (0,\pi/12)$, then (see (\ref{2.1}), (\ref{2.2})) 
\be \label{3.9} 
\begin{split}
& f_l(\alpha_0^{l,k_l,n})\in (0,1),\ l=1,\dots,9, \\ 
& f_l(U_n)\in (0,1),\ l=10,11,12. 
\end{split}
\ee  
Now, we make use the classical H. Hohr's theorem (1914) in the same way as it is done in our paper \cite{9}, Section 4, in order to generate of the following infinite sets 
\be \label{3.10} 
W_l^n,\ l=1,\dots,12,\ n=1,\dots,n_0;\ W_l^n\subset S_{\sigma_0^l,\delta}^l 
\ee  
(comp. (\ref{3.1})) of the elements 
\be \label{3.11} 
w_l^n\in W_l^n 
\ee  
as follows (see (\ref{1.1}), comp. \cite{9}) 
\be \label{3.12}  
\begin{split}
& |\zeta(w_l^n)|=f_l(\alpha_0^{l,k_l,n}),\ l=1,\dots,9, \\ 
& |\zeta(w_l^n)|=f_l(U_n),\ l=10,11,12, \\ 
& n=1,\dots,n_0, 
\end{split}
\ee  
where 
\be \label{3.13} 
\begin{split}
& w_l^n=w_l(U_n,\pi L,k_l,f_l,\sigma_0^l-\delta,\sigma_0^l+\delta),\ l=1,\dots,9, \\
& w_l^n=w_l(U_n,\pi L,f_l,\sigma_0^l-\delta,\sigma_0^l+\delta),\ l=10,11,12, 
\end{split}
\ee  
for every fixed set of admissible parameters. 

\begin{remark}
Corresponding set of grafts is defined by the equalities (\ref{3.12}). Of course, we choose only one element $w_l^n$ from every infinite set to construct the set of grafts (\ref{3.12}).  
\end{remark} 

\section{The set of elementary meta-functional equations} 

\subsection{} 

Now we make use the set of grafts (\ref{3.12}) (see also (\ref{2.1}), (\ref{2.2})) for grafting of every one from the $\zeta$-factorization formulas (\ref{2.3}) -- (\ref{2.11}), i.e. we use the substitutions (\ref{3.12}) in (\ref{2.3}) -- (\ref{2.11}). This procedure gives the following. 

\begin{mydef1}
For every fixed and admissible set of parameters 
\bdis 
\begin{split}
& U_n, \pi L, k_l, f_l, \sigma_0^l, \delta, \\ 
& n=1,\dots,n_0,\ l=1,\dots,12,\ L\geq L_0,\ 1\leq k_l\leq k_0 
\end{split}
\edis 
there is the set of following elementary meta-functional equations (exact forms): 
\be \label{4.1} 
|\zeta(w_1^n)|\prod_{r=1}^{k_1}\frac{\tilde{Z}^2(\alpha_r^{1,k_1,n})}{\tilde{Z}^2(\beta_r^{k_1,n})}=\frac 12-\frac 12|\zeta(w_{10}^n)|, 
\ee 

\be \label{4.2} 
|\zeta(w_2^n)|\prod_{r=1}^{k_2}\frac{\tilde{Z}^2(\alpha_r^{2,k_2,n})}{\tilde{Z}^2(\beta_r^{k_2,n})}=\frac 12+\frac 12|\zeta(w_{10}^n)|, 
\ee 

\be \label{4.3} 
|\zeta(w_3^n)|\prod_{r=1}^{k_3}\frac{\tilde{Z}^2(\alpha_r^{3,k_3,n})}{\tilde{Z}^2(\beta_r^{k_3,n})}=\frac 38-\frac 12|\zeta(w_{10}^n)|+\frac 18|\zeta(w_{11}^n)|, 
\ee  

\be \label{4.4} 
|\zeta(w_4^n)|\prod_{r=1}^{k_4}\frac{\tilde{Z}^2(\alpha_r^{4,k_4,n})}{\tilde{Z}^2(\beta_r^{k_4,n})}=\frac 38+\frac 12|\zeta(w_{10}^n)|+\frac 18|\zeta(w_{11}^n)|, 
\ee  

\be \label{4.5} 
|\zeta(w_5^n)|\prod_{r=1}^{k_5}\frac{\tilde{Z}^2(\alpha_r^{5,k_5,n})}{\tilde{Z}^2(\beta_r^{k_5,n})}=\frac{5}{16}-\frac{15}{32}|\zeta(w_{10}^n)|+\frac{3}{16}|\zeta(w_{11}^n)|-\frac{1}{32}|\zeta(w_{12}^n)|, 
\ee

\be \label{4.6} 
|\zeta(w_6^n)|\prod_{r=1}^{k_6}\frac{\tilde{Z}^2(\alpha_r^{6,k_6,n})}{\tilde{Z}^2(\beta_r^{k_6,n})}=\frac{5}{16}+\frac{15}{32}|\zeta(w_{10}^n)|+\frac{3}{16}|\zeta(w_{11}^n)|+\frac{1}{32}|\zeta(w_{12}^n)|, 
\ee 

\be \label{4.7} 
|\zeta(w_7^n)|\prod_{r=1}^{k_7}\frac{\tilde{Z}^2(\alpha_r^{7,k_7,n})}{\tilde{Z}^2(\beta_r^{k_7,n})}=
|\zeta(w_{10}^n)|, 
\ee  

\be \label{4.8} 
|\zeta(w_8^n)|\prod_{r=1}^{k_8}\frac{\tilde{Z}^2(\alpha_r^{8,k_8,n})}{\tilde{Z}^2(\beta_r^{k_8,n})}=
|\zeta(w_{11}^n)|, 
\ee  

\be \label{4.9} 
|\zeta(w_9^n)|\prod_{r=1}^{k_9}\frac{\tilde{Z}^2(\alpha_r^{9,k_9,n})}{\tilde{Z}^2(\beta_r^{k_9,n})}=
|\zeta(w_{12}^n)|.  
\ee 
\end{mydef1} 

\subsection{} 

By making use of (\ref{2.13}) (comp. \cite{8}, subsection 8.2) we obtain from our Theorem the following asymptotic forms of meta-functional equations (\ref{4.1}) -- (\ref{4.9}). 

\begin{mydef41}
\be \label{4.10} 
|\zeta(w_1^n)|\prod_{r=1}^{k_1}\frac{|\zeta(\frac 12+i\alpha_r^{1,k_1,n})|^2}{|\zeta(\frac 12+i\beta_r^{k_1,n})|^2}\sim \frac 12-\frac 12|\zeta(w_{10}^n)|, 
\ee  

\be \label{4.11} 
|\zeta(w_2^n)|\prod_{r=1}^{k_2}\frac{|\zeta(\frac 12+i\alpha_r^{2,k_2,n})|^2}{|\zeta(\frac 12+i\beta_r^{k_2,n})|^2}\sim \frac 12+\frac 12|\zeta(w_{10}^n)|, 
\ee 

\be \label{4.12} 
|\zeta(w_3^n)|\prod_{r=1}^{k_3}\frac{|\zeta(\frac 12+i\alpha_r^{3,k_3,n})|^2}{|\zeta(\frac 12+i\beta_r^{k_3,n})|^2}\sim \frac 38-\frac 12|\zeta(w_{10}^n)|+\frac 18|\zeta(w_{11}^n)|, 
\ee 

\be \label{4.13} 
|\zeta(w_4^n)|\prod_{r=1}^{k_4}\frac{|\zeta(\frac 12+i\alpha_r^{4,k_4,n})|^2}{|\zeta(\frac 12+i\beta_r^{k_4,n})|^2}\sim \frac 38+\frac 12|\zeta(w_{10}^n)|+\frac 18|\zeta(w_{11}^n)|, 
\ee  

\be \label{4.14} 
|\zeta(w_5^n)|\prod_{r=1}^{k_5}\frac{|\zeta(\frac 12+i\alpha_r^{5,k_5,n})|^2}{|\zeta(\frac 12+i\beta_r^{k_5,n})|^2}\sim \frac{5}{16}-\frac{15}{32}|\zeta(w_{10}^n)|+\frac{3}{16}|\zeta(w_{11}^n)|-\frac{1}{32}|\zeta(w_{12}^n)|, 
\ee  

\be \label{4.15} 
|\zeta(w_6^n)|\prod_{r=1}^{k_6}\frac{|\zeta(\frac 12+i\alpha_r^{6,k_6,n})|^2}{|\zeta(\frac 12+i\beta_r^{k_6,n})|^2}\sim \frac{5}{16}+\frac{15}{32}|\zeta(w_{10}^n)|+\frac{3}{16}|\zeta(w_{11}^n)|+\frac{1}{32}|\zeta(w_{12}^n)|, 
\ee  

\be \label{4.16} 
|\zeta(w_7^n)|\prod_{r=1}^{k_7}\frac{|\zeta(\frac 12+i\alpha_r^{7,k_7,n})|^2}{|\zeta(\frac 12+i\beta_r^{k_7,n})|^2}\sim |\zeta(w_{10}^n)|, 
\ee 

\be \label{4.17} 
|\zeta(w_8^n)|\prod_{r=1}^{k_8}\frac{|\zeta(\frac 12+i\alpha_r^{8,k_8,n})|^2}{|\zeta(\frac 12+i\beta_r^{k_8,n})|^2}\sim |\zeta(w_{11}^n)|, 
\ee 

\be \label{4.18} 
|\zeta(w_9^n)|\prod_{r=1}^{k_9}\frac{|\zeta(\frac 12+i\alpha_r^{9,k_9,n})|^2}{|\zeta(\frac 12+i\beta_r^{k_9,n})|^2}\sim |\zeta(w_{12}^n)|, 
\ee 
for $L\to\infty$. 
\end{mydef41} 

\subsection{} 

Now, we give some remarks. 

\begin{remark}
Let us notice the following: 
\begin{itemize}
	\item[(a)] formulae (\ref{4.1}) -- (\ref{4.9}) as well as formulae (\ref{4.10}) -- (\ref{4.18}) are synergetic ones, see our interpretation given in \cite{8}, 
	\item[(b)] the set of results of elementary interactions of the Riemann's zeta-function 
	\be \label{4.19} 
	\zeta(s), s=\sigma+it,\ \frac 12\leq \sigma<1 
	\ee 
	with itself is defined by the set of formulae (\ref{4.10}) -- (\ref{4.18}) (with respect to the sets of functions (\ref{2.1}) and (\ref{2.2})); that means results of interactions between the set 
	\be \label{4.20} 
	\left\{\left|\zeta\left(\frac 12+it\right)\right|^2\right\},\ t\leq L_0 
	\ee 
	and certain number of the following 12 sets 
	\be \label{4.21} 
	\{|\zeta(s)|\},\ s\in S_{\sigma_0^l,\delta}^l,\ l=1,\dots,12, 
	\ee 
	\item[(c)] (\ref{4.15}) gives, for example, the result of interactions of basic function (\ref{4.19}) with itself, namely, the result of the type 
	\bdis 
	\left\{\left|\zeta\left(\frac 12+it\right)\right|^2\right\},\ t\leq L_0,\ \{|\zeta(s)|\},\ s\in S_{\sigma_0^l,\delta}^l,\ l=9,10,11,12. 
	\edis 
\end{itemize} 
Let us call the set of these interactions as the $\zeta$-chemical reactions of the substances (\ref{4.20}) and (\ref{4.21}). This is then some analogue to what we have based on the classical Belousov-Zhabotinski chemical reaction, see \cite{8}. 
\end{remark} 

\begin{remark}
The set of elementary meta-functional equations (\ref{4.10}) -- (\ref{4.18}), as well as (\ref{4.1}) -- (\ref{4.9}), is the result (=set of $\zeta$-chemical compounds) of the corresponding nine types of $\zeta$-chemical reactions on the set of 13 substances (\ref{4.20}), (\ref{4.21}). 
\end{remark} 

\section{Examples of more complicated meta-functional equations} 

\subsection{} 

Here, we demonstrate the crossbreeding, see \cite{6} -- \cite{8}, on the subset (\ref{4.3}) -- (\ref{4.6}) of meta-functional equations as an example. Crossbreeding in this case means making use of Gauss' type elementary operations on this subset. namely, the first stage of the crossbreeding gives the following formulae: 
\be \label{5.1} 
\begin{split}
& (4.3) + (4.4)\Rightarrow \\ 
& |\zeta(w_3^n)||\prod_{r=1}^{k_3}\frac{\tilde{Z}^2(\alpha_r^{3,k_3,n})}{\tilde{Z}^2(\beta_r^{k_3,n})}+
|\zeta(w_4^n)||\prod_{r=1}^{k_4}\frac{\tilde{Z}^2(\alpha_r^{4,k_4,n})}{\tilde{Z}^2(\beta_r^{k_4,n})}=
\frac{3}{4}+\frac{1}{4}|\zeta(w_{11}^n)|, 
\end{split}
\ee 

\be \label{5.2} 
\begin{split}
	& (4.5) + (4.6)\Rightarrow \\ 
	& |\zeta(w_5^n)||\prod_{r=1}^{k_5}\frac{\tilde{Z}^2(\alpha_r^{5,k_5,n})}{\tilde{Z}^2(\beta_r^{k_5,n})}+
	|\zeta(w_6^n)||\prod_{r=1}^{k_6}\frac{\tilde{Z}^2(\alpha_r^{6,k_6,n})}{\tilde{Z}^2(\beta_r^{k_6,n})}=
	\frac{5}{8}+\frac{3}{8}|\zeta(w_{11}^n)|, 
\end{split}
\ee 
and elimination of $|\zeta(w_{11}^n)|$ from (\ref{5.1}) and (\ref{5.2}) represents the second stage. This gives the following. 

\begin{mydef42}
\be \label{5.3} 
\begin{split}
& 2\left\{
|\zeta(w_5^n)||\prod_{r=1}^{k_5}\frac{\tilde{Z}^2(\alpha_r^{5,k_5,n})}{\tilde{Z}^2(\beta_r^{k_5,n})}+
|\zeta(w_6^n)||\prod_{r=1}^{k_6}\frac{\tilde{Z}^2(\alpha_r^{6,k_6,n})}{\tilde{Z}^2(\beta_r^{k_6,n})}
\right\}+1=\\ 
& 3\left\{
|\zeta(w_3^n)||\prod_{r=1}^{k_3}\frac{\tilde{Z}^2(\alpha_r^{3,k_3,n})}{\tilde{Z}^2(\beta_r^{k_3,n})}+
|\zeta(w_4^n)||\prod_{r=1}^{k_4}\frac{\tilde{Z}^2(\alpha_r^{4,k_4,n})}{\tilde{Z}^2(\beta_r^{k_4,n})}
\right\}. 
\end{split}
\ee 
\end{mydef42} 

\subsection{} 

Secondary crossbreeding of the meta-functional equation (\ref{5.3}), for example, is the full analogue of the secondary crossbreeding on the class of corresponding complete hybrid formulas, see \cite{8}. First, (4.1) + (4.2) gives in the case 
\bdis 
k_1=k_2=k,\ 1\leq k\leq k_0 
\edis 
the result 
\be \label{5.4} 
|\zeta(w_1^n)|\prod_{r=1}^k\tilde{Z}^2(\alpha_r^{1,k,n})+|\zeta(w_2^n)|\prod_{r=1}^k\tilde{Z}^2(\alpha_r^{2,k,n})=\prod_{r=1}^k\tilde{Z}^2(\beta_r^{k,n}). 
\ee 
Second, we substitute (\ref{5.4}) into (\ref{5.3}) in the cases $k=k_5,k_6,k_3,k_4$ and obtain the following. 

\begin{mydef43}
\be \label{5.5} 
\begin{split}
& 2\left\{
|\zeta(w_5^n)|
\frac
{\prod_{r=1}^{k_5}\tilde{Z}^2(\alpha_r^{5,k_5,n})}
{|\zeta(w_1^n)|\prod_{r=1}^{k_5}\tilde{Z}^2(\alpha_r^{1,k_5,n})+
|\zeta(w_2^n)|\prod_{r=1}^{k_5}\tilde{Z}^2(\alpha_r^{2,k_5,n})}+\right. \\ 
& \left. +
|\zeta(w_6^n)|
\frac
{\prod_{r=1}^{k_6}\tilde{Z}^2(\alpha_r^{6,k_6,n})}
{|\zeta(w_1^n)|\prod_{r=1}^{k_6}\tilde{Z}^2(\alpha_r^{1,k_6,n})+
	|\zeta(w_2^n)|\prod_{r=1}^{k_6}\tilde{Z}^2(\alpha_r^{2,k_6,n})}
\right\}+1= \\ 
& 3\left\{
|\zeta(w_3^n)|
\frac
{\prod_{r=1}^{k_3}\tilde{Z}^2(\alpha_r^{3,k_3,n})}
{|\zeta(w_1^n)|\prod_{r=1}^{k_3}\tilde{Z}^2(\alpha_r^{1,k_3,n})+
	|\zeta(w_2^n)|\prod_{r=1}^{k_3}\tilde{Z}^2(\alpha_r^{2,k_3,n})}+\right. \\ 
& \left. +
|\zeta(w_4^n)|
\frac
{\prod_{r=1}^{k_4}\tilde{Z}^2(\alpha_r^{4,k_4,n})}
{|\zeta(w_1^n)|\prod_{r=1}^{k_4}\tilde{Z}^2(\alpha_r^{1,k_4,n})+
	|\zeta(w_2^n)|\prod_{r=1}^{k_4}\tilde{Z}^2(\alpha_r^{2,k_4,n})}
\right\}. 
\end{split}
\ee 
\end{mydef43} 

\begin{remark}
Of course, we obtain corresponding asymptotic formulae for (\ref{5.3}), (\ref{5.5}) if we use asymptotic meta-functional equations (\ref{4.10}) -- (\ref{4.15}) instead of (\ref{4.1}) -- (\ref{4.6}). 
\end{remark} 

\begin{remark}
Now we see that (\ref{1.2})=(\ref{4.15}) while (\ref{1.3}) is the asymptotic form of (\ref{5.3}). 
\end{remark}

I would like to thank Michal Demetrian for his moral support of my study of Jacob's ladders.

\end{document}